%rps.tex:
%%a Plain TeX file by Shalosh B. Ekhad and Doron Zeilberger (x pages)

%begin macros

\baselineskip=14pt
\parskip=10pt

\font\eightrm=cmr8 
\font\eighttt=cmtt8
\magnification=\magstephalf

\def\1{{\overline{1}}}
\def\2{{\overline{2}}}
\parindent=0pt
\overfullrule=0in

\def\frac#1#2{{#1 \over #2}}
%\headline={\rm  \ifodd\pageno  \RightHead  \else  \LeftHead  \fi}
%\def\RightHead{\centerline{
%Title
%}}
%\def\LeftHead{ \centerline{Doron Zeilberger}}
%end macros
\bf
\centerline
{
Automatic Solution of Richard Stanley's Amer. Math. Monthly Problem \#11610
}
\centerline
{
and ANY Problem of That Type
}
\rm
\bigskip
\centerline{ {\it
Shalosh B. EKHAD and
Doron 
ZEILBERGER}\footnote{$^1$}
{\eightrm  \raggedright
Department of Mathematics, Rutgers University (New Brunswick),
Hill Center-Busch Campus, 110 Frelinghuysen Rd., Piscataway,
NJ 08854-8019, USA.
%\break
{\eighttt zeilberg  at math dot rutgers dot edu} ,
\hfill \break
{\eighttt http://www.math.rutgers.edu/\~{}zeilberg/} .
Dec. 27, 2011.
Accompanied by Maple packages {\eighttt RPS} and {\eighttt RPSplus}
downloadable from Zeilberger's website.
Supported in part by the USA National Science Foundation.
}
}

{\bf Preamble}

Suppose you toss a (fair) coin $n$ times. If $n$ is large, the  {\it law of large numbers} promises you that
(with high probability)  you would {\it roughly} get as many Heads as Tails. But what is the {\it exact} probability that
you would have {\bf exactly} as many Heads as Tails? If $n$ is odd, the answer is easy ({\bf you do it!}).
If $n$ is even, then it is almost as easy, and there is a {\it nice}, ``closed-form'' formula for that
probability, namely $n!/((n/2)!^2 2^n)$.

Richard Stanley [St1] proposed the problem of finding, $a(n)$,
the number of $n$-letter words in  the alphabet $\{H,T\}$ where
there are as many occurrences of ``HT'' (i.e. Head immediately followed by Tail) as there are occurrences
of ``TT'' (two Tails in a row). He didn't give a ``closed form'' formula, but he gave something almost as good,
an {\it explicit} formula as an ({\it algebraic}, as it turned out, in fact quadratic)
formal power series for the (ordinary) generating function $P(t):=\sum_{n=0}^{\infty} a(n)t^n$.

The fact that the generating function, $P(t)$, is an {\it algebraic} generating function is not at all surprising! 
This can be seen in (at least) two ways.

One way is to show that the ``language'' of words with as many occurrences of ``HT'' as
``HH'' is {\it context-free} (type 2) with an {\it unambiguous grammar}, and hence its {\it weight-enumerator}
is {\it algebraic}. It is possible to ({\it automatically}!) generate its grammar, and then {\it automatically}
generate a system of algebraic equations one of whose unknowns is the desired generating function,
and solving that system would (presumbly, we didn't do it) yield Stanley's proposed expression.

A better way is to find (automatically, of course!), the {\it rational generating function}
$F(t;z[HT],z[TT])$ that is the {\it weight-enumerator} of all words in the alphabet $\{H,T \}$ according
to the weight $Weight(w)=t^{length(w)}z[HT]^{\#HT(w)}z[TT]^{\#TT(w)}$. This can be done in several ways,
including the {\it Goulden-Jackson method}, beautifully surveyed in [NZ], and efficiently implemented
in the Maple package

{\tt http://www.math.rutgers.edu/\~{}zeilberg/tokhniot/DAVID\_IAN}

accompanying that article.

Having done that, the desired generating function, $P(t)$, is the coefficient of $s^0$ (i.e. {\it the constant term}) in
$F(t;s,1/s)$.
Hillel Furstenberg[F] promises us that $P(t)$ is an algebraic formal power series in $t$,
and his proof implies a (rather awkward and inefficient) algorithm (using Cauchy's integral formula and residues) for computing it.
This should yield a second {\it rigorous} derivation of Stanley's proposed solution.

But since we know {\it a priori} (by ``general nonsense'') that the desired sequence 
belongs to the {\it algebraic ansatz} (see [Z1] and the {\bf wonderful} new book by Manuel Kauers and
Peter Paule[KP], that should be {\it required reading} to {\it any} mathematics student (and professional!))
a {\it semi-rigorous} derivation would be to crank out the first $40$ (or even fewer) terms
of Stanley's sequence (a quick way would be to find the first $40$ terms in the expansion of the
constant term, in $s$,
of $F(t;s,1/s)$ above), and then use a {\it guessing} program, e.g. {\tt listoalg}
in the Maple package {\it gfun}([SaZ]) (but {\it please} enlarge the very small default values of the parameters)
that now is part of Maple, or procedure {\tt Empir} in my own Maple package 
{\tt http://www.math.rutgers.edu/\~{}zeilberg/tokhniot/SCHUTZENBERGER} .
 
In order to make the above {\it semi-rigorous} derivation {\it fully rigorous} (for those obtuse people
who desire it), one would need to derive {\it a priori} bounds on the degree (in $t$ and $P(t)$) of the defining
equation $F(t,P(t))=0$ for the desired generating function $P(t)$. Unlike the $C$-finite ansatz
(see [Z2] and [KP]) where finding these upper bounds is trivial, I don't know how to do it in the present case.
But there is {\bf another} way to make everything fully rigorous. Via the {\it holonomic ansatz}!

Using the {\it Continuous Almkvist-Zeilberger Algorithm}[AlZ], that is implemented in procedure {\tt AZc} of
the Maple package {\tt http://www.math.rutgers.edu/\~{}zeilberg/tokhniot/EKHAD}, one can obtain a
{\it differential equation} (and its {\it proof} (a certain {\it certificate})), and then
verify that the above ``conjectured'' algebraic expression for $P(t)$ satisfies that very same differential equation,
and check that the  initial conditions match.

{\bf The general case}

The beauty of {\it algorithmic} mathematics is that it is not much harder to write a {\it general} program to handle
a whole class of problems rather than just solve {\it one} problem. The above discussion applies equally
to any (finite) {\it alphabet} (not just a two-lettered one) and any two {\it distinguished} substrings,
$w_1$ and $w_2$ not just $HT$ and $TT$.

{\bf The Maple package RPS}

Since we require procedures from four different Maple packages ({\tt DAVID\_IAN}, {\tt SCHUTZENBERGER},
{\tt EKHAD} and {\tt AsyRec}), we conveniently assembled all the necessary procedures,
together with new ``interfacing code'' needed to solve problems of the above type.
The result is the Maple package {\tt RPS}
(named after Richard Peter Stanley), available free of charge from:

{\tt http://www.math.rutgers.edu/\~{}zeilberg/tokhniot/RPS} .

This Maple package does much more! It computes {\it holonomic representations}
(or as Richard Stanley [St2] would say, {\it P-recursive} ones), that are used, in turn, to
derive {\it asymptotic expressions} using procedures borrowed from

{\tt http://www.math.rutgers.edu/\~{}zeilberg/tokhniot/AsyRec} \quad .

Hence we have ``three-quarters'' of the Kauers-Paule ``concrete tetrahedron'': generating function, recurrence,
and asymptotics. The last one ``definite sum'' could also be obtained, but we would (usually) get
complicated and ugly {\it multi-sums} with many sigmas, so it would be stupid to look for these.

Out of sheer laziness we have only programmed the case where the two distinguished words, $w_1$, $w_2$, have the same length.

{\bf Precomputed Output of the Maple package RPS}

The ``front'' of this article, the webpage

{\tt http://www.math.rutgers.edu/\~{}zeilberg/mamarim/mamarimhtml/rps.html} \quad ,

contains links to several ``webbooks'' that systematically states (proved!) (algebraic) generating functions,
recurrence equations, and asymptotics for analogs of Stanley's problem for
{\it all} possible pairs of words (up to trivial images under permutations of the letters)
of the same length (let's call it $k$) for an $m$-letter alphabet for the following cases.

$\bullet$ $m=2,k=2$: {\tt http://www.math.rutgers.edu/\~{}zeilberg/tokhniot/oRPS22} (containing $3$ propositions)

$\bullet$ $m=2,k=3$: {\tt http://www.math.rutgers.edu/\~{}zeilberg/tokhniot/oRPS23} (containing $11$ propositions)

$\bullet$ $m=2,k=4$: {\tt http://www.math.rutgers.edu/\~{}zeilberg/tokhniot/oRPS24} (containing $38$ propositions)

$\bullet$ $m=3,k=2$: {\tt http://www.math.rutgers.edu/\~{}zeilberg/tokhniot/oRPS32} (containing $6$ propositions)

$\bullet$ $m=3,k=3$: {\tt http://www.math.rutgers.edu/\~{}zeilberg/tokhniot/oRPS32} (containing $40$ propositions)

$\bullet$ $m=4,k=2$: {\tt http://www.math.rutgers.edu/\~{}zeilberg/tokhniot/oRPS42} (containing $7$ propositions)

$\bullet$ $m=4,k=3$: {\tt http://www.math.rutgers.edu/\~{}zeilberg/tokhniot/oRPS43} (containing $63$ propositions)

$\bullet$ $m=5,k=2$: {\tt http://www.math.rutgers.edu/\~{}zeilberg/tokhniot/oRPS52} (containing $7$ propositions)

$\bullet$ $m=5,k=3$: {\tt http://www.math.rutgers.edu/\~{}zeilberg/tokhniot/oRPS53} (containing $69$ propositions)

$\bullet$ $m=6,k=2$: {\tt http://www.math.rutgers.edu/\~{}zeilberg/tokhniot/oRPS62} (containing $7$ propositions)

$\bullet$ $m=6,k=3$: {\tt http://www.math.rutgers.edu/\~{}zeilberg/tokhniot/oRPS63} (containing $70$ propositions)

For Neil Sloane's sake, we have also computed the first $50$ terms of each of the considered enumerating sequences.
All the sequences for $m=2$ and $k=2,3,4$ have already been  entered to OEIS by R.H. Hardin, for example
{\tt http://oeis.org/A164147}. Some of the pages for these sequences come with conjectured recurrences. The present
webbooks supply rigorous proofs to all them, and supplies proved recurrences for the remaining ones.

{\bf The Maple Package RPSplus}

With hardly any more (programming) effort, one can consider the enumerating
sequences of words for which, for three given positive integers $a_1$,$a_2$ and $r$,

``$a_1$ times [the numer of occurrences of $w_1$]''  {\bf minus} $a_2$ times ``[the numer of occurrences of $w_2$] ``
{\bf equals} $r$.

Once again the generating functions are
guaranteed to be algebraic and everything goes through. See the procedures listed in
{\tt ezraG();} in the more general Maple package {\tt RPSplus}, available from

{\tt http://www.math.rutgers.edu/\~{}zeilberg/tokhniot/RPSplus } \quad .

Readers are welcome to generate their own output.

{\bf What About Several Distinguished Words?}

The more general case where one has a finite alphabet of, say, $m$ letters, and
$s$, say, distinguished words $w_1, \dots , w_s$, and $v$ diophantine affine linear
relations between the quantities ``number of occurrences of $w_i$'', then we leave the {\it algebraic ansatz}
and enter the {\it holonomic} ansatz. By WZ theory we are {\it guaranteed} that the
enumerating sequence, in each case, is holonomic (alias $P$-recursive), and we are justified,
semi-rigorously, just to {\it guess} the holonomic description, using
{\tt gfun}'s {\tt listtorec}, or procedure {\tt Findrec} in the Maple package {\tt RPS}.

For those obtuse people who insist on a {\it rigorous} proof, they are welcome to use
the Maple package

{\tt http://www.math.rutgers.edu/\~{}zeilberg/tokhniot/MultiAlmkvistZeilberger} \quad ,

that is one of the Maple packages accompanying the seminal article [ApZ].
Alas, it may take them quite some time, and frankly, for us, a {\it semi-rigorous} proof suffices.
But so far, we ran out of steam, and we do not even have an implementation of
the semi-rigorous, pure guessing, version.

\vfill\eject

{\bf References}

[AlZ] Gert Almkvist and Doron Zeilberger, {\it The Method of Differentiating Under The Integral Sign},
J. Symbolic Computation {\bf 10}(1990), 571-591.

[ApZ] Moa Apagodu and Doron Zeilberger,
{\it Multi-Variable Zeilberger and Almkvist-Zeilberger Algorithms and the Sharpening of Wilf-Zeilberger Theory},
Adv. Appl. Math. {\bf 37}(2006), (Special issue in honor of Amitai Regev), 139-152.

[F] Hillel Furstenberg, {\it Algebraic Functions over finite fields}, J. Algebra {\bf 7} (1967). 271-277.

[KP] Manuel Kauers and Peter Paule, {\it ``The Concrete Tetrahedron''}, Springer, 2011.

[NZ] John Noonan and Doron Zeilberger, {\it The Goulden-Jackson Cluster Method: Extensions, Applications, and Implementations},
J. Difference Eq. Appl. {\bf 5}(1999), 355-377.
\hfill\break 
{\tt http://www.math.rutgers.edu/\~{}zeilberg/mamarim/mamarimhtml/gj.html}

[SaZ] Bruno Salvy and Paul Zimmermann,
{\it Gfun: a Maple package for the manipulation of generating and holonomic functions in one variable},
ACM Transactions on Mathematical Software, {\bf 20} (1994), 163-177.

[St1] Richard P. Stanley, {\it Problem  \#11610}, Amer. Math. Monthly {\bf 118(10)} (Dec. 2011), 937.

[St2] Richard P. Stanley, {\it Differentiably finite power series}, Europ. J. Combinatorics {\bf 1} (1980), 175-188.

[Z1] Doron Zeilberger, {\it An Enquiry Concerning Human (and Computer!) [Mathematical] Understanding},
in: ``RANDOMNESS AND COMPLEXITY, FROM LEIBNIZ TO CHAITIN'', Cristian C. Calude, ed.,
\hfill\break 
{\tt http://www.math.rutgers.edu/\~{}zeilberg/mamarim/mamarimhtml/enquiry.html}

[Z2] Doron Zeilberger, {\it The C-finite Ansatz}, submitted, \hfill\break
{\tt http://www.math.rutgers.edu/\~{}zeilberg/mamarim/mamarimhtml/cfinite.html} .

\end